\documentclass{article}
\usepackage{geometry}
\usepackage{amssymb}
\usepackage{amsmath}

\newtheorem{theorem}{Theorem}[section]

\newcommand{\R}{\mathbb{R}}

\newcommand{\C}{\mathbb{C}}

\renewcommand{\i}{\mathrm{i}}

\newcommand{\const}{\mathrm{const}}
\newcommand{\spn}{\operatorname{span}}

\newcommand{\la}{\lambda}

\title{On deformations of the dispersionless Hirota equation}
\author{
Wojciech Kry\'nski\thanks{
Institute of Mathematics, Polish Academy of Sciences, ul.~\'Sniadeckich 8, 00-959 Warszawa, Poland
\newline 
E-mail: krynski@impan.pl
}
}

\begin{document}
\maketitle
\begin{abstract}
The hyper-CR Einstein-Weyl structures on $\R^3$ can be described in terms of the solutions to the dispersionless Hirota equation. In the present paper we show that simple geometric constructions on the associated twistor space lead to deformations of the Hirota equation that have been introduced recently by B.~Kruglikov and A.~Panasyuk. Our method produces also the hyper-CR equation and can be applied to other geometric structures related to different twistor constructions.
\end{abstract}

\section{Introduction}
It is proved in \cite{DK} that solutions $f\colon \R^3\to\R$ to the dispersionless Hirota equation
\begin{equation}\label{eq1}
af_1f_{23}+bf_2f_{13}+cf_3f_{12}=0,
\end{equation}
where $a+b+c=0$, are in a correspondence with the hyper-CR Einstein-Weyl structures on $\R^3$. In the recent paper \cite{KP} it is shown that the hyper-CR class can be equivalently described by four other families (A)-(D) of integrable dispersionless PDEs that are only B\"acklund equivalent to \eqref{eq1}. The equations are as follows
\begin{align}
&(\la_2(p_2)-\la_3(p_3))f_1f_{23}+(\la_3(p_3)-\la_1(p_1))f_2f_{13} +(\la_1(p_1)-\la_2(p_2))f_3f_{12}=0,\tag{A}\\
&f_1f_{13}-f_3f_{11}+(\la_2(x_2)-\la_3(p_3))(f_1f_{23}-f_2f_{13})+\la_2'(p_2)f_1f_3=0,\tag{B}\\
&f_1f_{13}-f_3f_{11}+e^{-\la_3'(p_3)p_2}(f_2f_{12}-f_1f_{22})+\la_3''(p_3)p_2f_1^2=0,\tag{C}\\
&(a(p_1,p_2)-\la_3(p_3))(f_1f_{23}-f_2f_{13})+b(p_1,p_2)(f_3f_{11}+f_3f_{22}-f_1f_{13}-f_2f_{23})=0,\tag{D}
\end{align}
where $p_i$, $i=1,2,3$, are coordinate functions on $\R^3$, the functions $\la_i$, $i=1,2,3$, depend on one coordinate function only and the functions $a$ and $b$ depend on two coordinate functions and satisfy $a_1=b_2$ and $a_2=-b_1$. The derivation of equations (A)-(D) is based on an analysis of the Nijenhuis tensors associated to the Veronese webs that underlay the hyper-CR structures. All these equations fit into the general scheme that includes also the dKP equation \cite{DMT, GN} and the Manakov-Santini system \cite{MS} and is covered by \cite{FK} where a characterization of 3-dimensional integrable dispersionless equations in terms of the Einstein-Weyl geometry is given (see also \cite{FM}).

A purpose of the present note is to explore ideas of \cite[Section 3]{DK} and interpret the results of \cite{KP} from the point of view of the geometry of the twistor space associated to the Einstein-Weyl structures. Moreover, we extend the results of \cite{KP} and show that the hyper-CR equation
\begin{equation}\label{eq4}
H_{13}-H_{22}+H_2H_{33}-H_3H_{23}=0
\end{equation}
considered in \cite{D} can be also derived from the Hirota equation by simple geometric constructions on the twistor space. The hyper-CR equation is missed in the approach of \cite{KP} since it is not directly related to any Nijenhuis tensor. Another advantage of our approach is that it immediately solves the problem of realizability which was crucial in \cite{KP}.

The constructions presented in the present paper can be also applied to other dispersionless equations related to different twistorial constructions. We analyse two examples at the end of the paper in Sections 5 and 6. Namely we consider systems describing the hyper-Hermitian structures and the Veronese webs in dimension 4 (which easily generalize to higher dimensions).

In the hyper-Hermitian case the systems descent to the Pleba\'nski equations provided that the structures are hyper-K\"ahler. In this case similar ideas can be found in \cite{DM} (see also \cite{C}).

The Veronese webs can be considered as a special class of $GL(2)$-structures which are of the utmost importance for the geometry of ordinary differential equations (see \cite{B,DT,FK1,GN1,K1,N} and references therein). It is shown in \cite{K1} that the Veronese webs are described by a hierarchy of integrable systems that generalizes the Hirota equation (integrable systems describing more general classes of $GL(2)$-structures are given in \cite{FK1,KM}). In the present paper we prove that the webs can be described by two other hierarchies. The systems are defined as the compatibility conditions for the system
\begin{equation}\label{sys0}
f_i+H_{i+1}f_k=0, \qquad i=0,\ldots,k-1.
\end{equation}
The two hierarchies generalize the most symmetric equation in the family (C) on the one hand and the hyper-CR equation on the other hand.

We deal with the most symmetric equations in the families (A)-(C) in the next section. They correspond to constant functions $\la_i$. The deformations with arbitrary $\la_i$ are presented in Section 3. The family (D) is treated separately in Section 4 since it is related to the complex twistor space whereas the families (A)-(C) are related to the real twistor space. 

\section{Symmetric equations in families (A), (B) and (C)}
Let $T$ be the real twistor space of a hyper-CR Einstein-Weyl structure on a manifold $M$ and let $\pi\colon T\to\R P^1$ be the corresponding fibration over the projective space (c.f. \cite{D,DK}). Points in $T$ are surfaces in $M$ and each fiber of $\pi$ can be interpret as a foliation of $M$. There is a one-parameter family of the foliations parameterized by $\R P^1$ and the family appears to be a Veronese web (c.f. \cite{GZ,GZ1,Z}) as explained in \cite{DK}.  The leaves of the foliations are exactly the totally-geodesic surfaces of the Einstein-Weyl structure. Note that the twistor space $T$ is also the twistor space of the Veronese web considered in \cite{GZ1,Z}.

Pick three points $\la_1,\la_2,\la_3\in \R P^1$ and take a curve $\gamma_p\colon \R P^1\to T$ corresponding to a point $p$ in $M$. We define coordinates on $M$ by the following formula
$$
p=(p_1,p_2,p_3)=(\gamma_p(\la_1),\gamma_p(\la_2),\gamma_p(\la_3)),
$$
where $\gamma_p(\la_1)$, $\gamma_p(\la_2)$, $\gamma_p(\la_3)$ are points of the intersection of $\gamma_p$ with the three fibers $\pi^{-1}(\la_1)$, $\pi^{-1}(\la_2)$, $\pi^{-1}(\la_3)$. We assume here that the fibres are identified with $\R$ in some arbitrary way, thus we treat $p_i$, $i=1,2,3$, as functions on $\R^3$ taking values in $\R$. 

Note that foliations $p_i=\const$ correspond exactly to the fibers $\pi^{-1}(\lambda_i)$, $i=1,2,3$. Now, let $\la_4\in \R P^1$ be a fourth point in $\R P^1$ and let $f$ be a function on $M$ such that $f=\const$ defines the foliation $\pi^{-1}(\la_4)$, i.e.
$$
f(p)=\gamma_p(\la_4).
$$

For simplicity we assume that all $\la_i$, $i=1,2,3$, are in the real line $\R\subset\R P^1=\R\cup\{\infty\}$. If additionally $\la_4\neq\infty$ then the tangent bundle of the upper mentioned Veronese web that underlies the hyper-CR Einstein-Weyl is annihilated by the following $\la$-dependent one-form as explained in \cite{DK}
\begin{equation}\label{formV31}
\begin{aligned}
\alpha_\la&=(\la_4-\la_1)(\la-\la_2)(\la-\la_3)f_1dp_1+(\la-\la_1)(\la_4-\la_2)(\la-\la_3)f_2dp_2\\
&+(\la-\la_1)(\la-\la_2)(\la_4-\la_3)f_3dp_3,
\end{aligned}
\end{equation}
where $f_i=\frac{\partial f}{\partial p_i}$. If $\la_4=\infty$ then the one-form takes the form
\begin{equation}\label{formV32}
\alpha_\la=(\la-\la_2)(\la-\la_3)f_1dp_1+(\la-\la_1)(\la-\la_3)f_2dp_2+(\la-\la_1)(\la-\la_2)f_3dp_3.
\end{equation}
In both cases the one-form $\alpha_\la$ satisfies the following integrability condition
$$
\alpha_\la\wedge d\alpha_\la=0
$$
which written in terms of $f$ is equivalent to the Hirota equation \eqref{eq1}, with
$$
a=(\la_4-\la_1)(\la_2-\la_3),\quad b=(\la_4-\la_2)(\la_3-\la_1),\quad c=(\la_4-\la_3)(\la_1-\la_2),
$$
provided that all $\la_i$, $i=1,2,3,4$, lay on the real line in $\R P^1$ or with
$$
a=(\la_2-\la_3),\quad b=(\la_3-\la_1),\quad c=(\la_1-\la_2),
$$
if $\la_4=\infty$.

Equation \eqref{eq1} is the most symmetric equation in the family (A). Above, we assumed that all points $\la_i$, $i=1,2,3,4$, are different.  Now, we would like to consider degenerate cases when some of them coincide. 

\paragraph{Case 1: $\la_1$ and $\la_2$ coincide.} Assume $\la_4=\infty$ and let 
$$
q_1= p_1,\quad q_2= \frac{p_2-p_1}{\la_2-\la_1},\quad q_3= p_3
$$
be new coordinartes on $M$. Then, denoting $\delta=\la_2-\la_1$, we can write
$$
a=\delta-b,\quad c=-\delta.
$$
Moreover
$$
p_1=q_1, \quad p_2=q_1+\delta q_2,\quad p_3=q_3,
$$
and in the new coordinates the Hirota equation takes the form
$$
\delta^{-2}(f_3f_{22}-f_2f_{23}+b(f_2f_{13}-f_1f_{23})) + \delta^{-1}(f_1f_{23}-f_3f_{12})=0.
$$
Multiplying by $\delta^2$ in the limit $\delta\to 0$ we get
\begin{equation}\label{eq2}
f_3f_{22}-f_2f_{23}+b(f_2f_{13}-f_1f_{23})=0.
\end{equation}
which is (up to a permutation of coordinates) the most symmetric equation in the family (B). Note that
$$
\lim_{\delta\to 0}q_2 = \gamma_p'(\la_1).
$$
Thus, if $\delta\to 0$ then the coordinates of a point $p\in M$ are $(\gamma_p(\la_1), \gamma_p'(\la_1), \gamma_p(\la_3))$.

\paragraph{Case 2: $\la_1$, $\la_2$ and $\la_3$ coincide.} Assume $\la_4=\infty$ and let 
$$
r_1= p_1,\quad r_2= \frac{p_2-p_1}{\la_2-\la_1},\quad r_3=\frac{p_3-2p_2+p_1}{(\la_2-\la_1)(\la_3-\la_2)}
$$
be new coordinates on $M$. Then, denoting $\delta=\la_2-\la_1$ and $\gamma=\la_3-\la_2$, we can write
$$
a=-\gamma,\quad b=\delta+\gamma,\quad c=-\delta
$$
and
$$
p_1=r_1,\quad p_2=r_1+\delta r_2,\quad p_3=r_1+2\delta r_2+\delta\gamma r_3.
$$
Hence, in the new coordinates the Hirota equation takes the form
$$
\begin{aligned}
&\frac{1}{\delta^2}(f_2f_{13}-f_1f_{23})+\frac{1}{\delta\gamma}(f_2f_{13}-f_3f_{12})\\
&+\frac{1}{\delta^3}f_2f_{23}+\frac{1}{\delta^2\gamma}(2f_1f_{33}+f_3f_{22}-2f_3f_{13}-f_2f_{23})-\frac{1}{\delta\gamma^2}f_2f_{23}\\
&-\frac{1}{\delta^3\gamma}(f_3f_{23}+f_2f_{33})+\frac{1}{\delta^2\gamma^2}(f_2f_{33}-f_3f_{23})+\frac{2}{\delta\gamma^3}f_3f_{23}=0
\end{aligned}
$$
Assuming that $\delta=\gamma$, multiplying by $\delta^3$ and rescaling $(r_1,r_2,r_3)\mapsto (r_1,r_2,2r_3)$ in the limit $\delta\to 0$ we get
\begin{equation}\label{eq3}
f_3f_{22}-f_2f_{23}+f_1f_{33}-f_3f_{13}=0
\end{equation}
which is (up to a permutation of coordinates) the most symmetric equation in the family (C). Note that
$$
\lim_{\delta=\gamma\to 0}r_2=\gamma_p'(\la_1),\qquad \lim_{\delta=\gamma\to 0}r_3=\gamma_p''(\la_1).
$$
Thus, if $\delta\to 0$ then the coordinates of a point $p\in M$ are $(\gamma_p(\la_1), \gamma_p'(\la_1), \gamma_p''(\la_1))$.

\paragraph{Case 3: $\la_1$, $\la_2$, $\la_3$ and $\la_4$ coincide.} We consider the coordinates $r=r(p)$ as in the previous case but instead of $f(p)=\gamma_p(\la_4)$ we take
$$
H(p)=\frac{1}{\delta^3}(f(p)-3p_3+3p_2-p_1)
$$
where, for simplicity, we assume $\delta=\la_2-\la_1=\la_3-\la_2=\la_4-\la_3$. Then
$$
a=-3\delta^2,\quad b=4\delta^2,\quad c=-\delta^2,
$$
and in the new coordinates
$$
\begin{aligned}
\frac{1}{\delta^7}(-3f_1f_{23}+4f_2f_{13})+\frac{1}{\delta^8}(6f_1f_{33}-f_2f_{23}-6f_3f_{13}+f_3f_{22})
+\frac{1}{\delta^9}(-2f_2f_{33}+2f_3f_{23})=0
\end{aligned}
$$
additionally, since $f=\delta^3H+3\delta^2q_3+3\delta q_2+q_1$
$$
f_1=\delta^3H_1+1,\qquad f_2=\delta^3H_2+3\delta,\qquad f_3=\delta^3H_3+3\delta^2.
$$
Thus, substituting $H$ instead of $f$, multiplying by $\delta^6$ and rescaling $(r_1,r_2,r_3,H)\mapsto (r_1,r_2,2r_3,6H)$ in the limit $\delta\to 0$ we get exactly the hyper-CR equation \eqref{eq4}. Note that
$$
\lim_{\delta\to 0}H=\gamma_p'''(\la_1).
$$

\paragraph{Remark.} 
Equations \eqref{eq4} and \eqref{eq3} are written down in the same coordinate system on $M$. They are equivalent in the following sense. Consider the system
$$
f_1+H_2f_3=0,\qquad f_2+H_3f_3=0.
$$
If $f$ is given then a solution $H$ exists if and only if $f$ solves \eqref{eq3}. On the other hand, if $H$ is given then a solution $f$ exists if and only if $H$ solves \eqref{eq4}.

\section{Deformations of families (A), (B) and (C)}
The Hirota equation can be deformed further in the following way. Instead of the fibres $\pi^{-1}(\la_1)$, $\pi^{-1}(\la_2)$ and $\pi^{-1}(\la_3)$ one can consider three one-dimensional submanifolds $P_1$, $P_2$ and $P_3$ of the twistor space $T$ that are transversal to curves $\gamma_p\colon \R P^1\to T$ corresponding to points in $M$. Precisely, we consider an open subset of $M$ consisting of points $p$ such that $\gamma_p$ intersects transversally $P_i$, $i=1,2,3$. The three submanifolds $P_1$, $P_2$ and $P_3$ define coordinates on $M$ by the formula
\begin{equation}\label{coordinates}
p=(p_1,p_2,p_3)=(\gamma_p\cap P_1, \gamma_p\cap P_2, \gamma_p\cap P_3)
\end{equation}
provided that a parametrization $P_i\simeq\R$ is given for all $i=1,2,3$. 

For a point $p\in M$ define
\begin{equation}\label{functions}
\la_1(p)=\pi(\gamma_p\cap P_1),\quad \la_2(p)=\pi(\gamma_p\cap P_2),\quad \la_3(p)=\pi(\gamma_p\cap P_3).
\end{equation}
Thus $\la_1(p),\la_2(p),\la_3(p)\in\R P^1$. We shall assume that there is an affine parameter on $\R P^1$ and $\la_i$, $i=1,2,3$, are away from $\infty\in\R P^1$. Moreover, we shall assume $\la_4=\infty$. We get that the Veronese web underlying the Einstein-Weyl structure is defined by the one-form $\alpha_\la$ given by \eqref{formV32}, where now $\la_i$, $i=1,2,3$, are functions on $M$ with values in $\R$. However, formulae \eqref{coordinates} and \eqref{functions} combined imply that each of $\la_i$, $i=1,2,3$, depends on one coordinate function only, precisely
$$
\la_1=\la_1(p_1),\quad \la_2=\la_2(p_2),\quad \la_3=\la_3(p_3).
$$
The integrability condition $\alpha_\la\wedge d\alpha_\la=0$ reads
\begin{equation}\label{eqd1}
(\la_2(p_2)-\la_3(p_3))f_1f_{23}+(\la_3(p_3)-\la_1(p_1))f_2f_{13}+(\la_1(p_1)-\la_2(p_2))f_3f_{12}=0
\end{equation}
This is exactly a general equation in the family (A).

Note that if $P_i$ projects regularly to $\R P^1$ then an affine parameter on $\R P^1$ defines a natural parameter on $P_i$. Later, without lost of generality, we shall assume that all $P_i$, $i=1,2,3$, are parameterized such that the corresponding functions $\la_i$ are linear in $p_1$, $p_2$ and $p_3$, respectively (compare \cite[Remark 4.2]{KP}). In the consecutive paragraphs we shall consider deformations of \eqref{eqd1}. Namely, we will assume that $P_1$ coincides with $P_2$ and $P_3$. In this way we will get the general equations from families (B) and (C).

\paragraph{Case 1: $P_1$ and $P_2$ coincide.}
We assume that both $P_1$ and $P_2$ are graphs of linear functions that are mutually parallel and differ by $\delta$, i.e. for a point $p\in M$ we have
$$
\la_1(p)=Ap_1+B,\quad \la_2(p)=Ap_2+B+\delta,
$$
for some $A,B\in\R$. Denote
$$
\tilde\delta=\la_2(p)-\la_1(p)
$$
and introduce new coordinates
$$
q_1=p_1, \quad q_2=\frac{p_2-p_1}{\tilde\delta}, \quad q_3=p_3.
$$
Thus
$$
\tilde\delta=\frac{\delta}{1-Aq_2}
$$
and
$$
p_1=q_1,\quad p_2=q_1+\tilde\delta q_2,\quad p_3=q_3.
$$
Note that $q_2$ approximates $\gamma_p'(\la_1(p))$, i.e. $\lim_{\delta\to 0}q_2=\gamma'_p(\la_1(p))$, where $\gamma_p$ is the curve in $T$ corresponding to $p\in M$. Using the chain rule we compute the following transformations of the derivatives of $f$ under the change of coordinates $p\mapsto q(p)$ 
$$
f_1\mapsto f_1-\psi^{-1}f_2,\quad f_2\mapsto\psi^{-1}f_2,\quad f_3\mapsto f_3
$$
and
$$
\begin{aligned}
&f_{12}\mapsto \psi^{-1}f_{12}-\psi^{-2}f_{22}+\psi_2\psi^{-2}f_2\\
&f_{13}\mapsto f_{13}-\psi^{-1}f_{23}\\
&f_{22}\mapsto\psi^{-2}f_{22}-\psi_2\psi^{-3}f_2\\
&f_{23}\mapsto\psi^{-1}f_{23} 
\end{aligned}
$$
where we denote
$$
\psi=\tilde\delta+q_2\frac{d\tilde\delta}{dq_2}=\frac{\delta}{(1-Aq_2)^2}.
$$
Hence, we get that in the new coordinates, after multiplying by $\delta^2$ and passing to the limit $\delta\to 0$, equation \eqref{eqd1} takes the form
\begin{equation}\label{eqd2}
(1-Aq_2)(f_3f_{22}-f_2f_{23})+(\la_3(q_3)-\la_1(q_1))(f_2f_{13}-f_1f_{23})-2Af_2f_3=0.
\end{equation}
Now, the transformation
$$
(1-Aq_2)\mapsto e^{-Aq_2}
$$
simplifies \eqref{eqd2} to
\begin{equation}\label{eqd2a}
f_3f_{22}-f_2f_{23}+(\la_3(q_3)-\la_1(q_1))(f_2f_{13}-f_1f_{23})-Af_2f_3=0.
\end{equation}
That is the most general equation in the family (B) of \cite{KP} (having in mind \cite[Remark 4.2]{KP}).

\paragraph{Case 2: $P_1$, $P_2$ and $P_3$ coincide.} We consider the coordinates $q=(q_1,q_2,q_3)$ from the previous case but assume now that the submanifold $P_3$ is parallel to $P_1$ and
$$
\la_1(q)=Aq_1+B,\quad \la_3(q)=Aq_3+B+\delta,
$$
for some $A,B\in\R$. Note that $q_1=p_1$ and $q_3=p_3$ are the original coordinates defined by \eqref{coordinates}. We introduce new coordinates
$$
r_1=q_1,\quad r_2=q_2,\quad r_3=2\frac{q_3-\tilde\delta q_2-q_1}{\tilde\delta^2}
$$
where
$$
\tilde\delta=\frac{\delta}{1-Aq_2}=\frac{\delta}{1-Ar_2}.
$$
Since $q_2=r_2=\gamma_p'(\la_1(p))$ it can be easily shown that $r_3$ approximates $\gamma''_p(\la_1(p))$, i.e. $\lim_{\delta\to 0}r_3=\gamma''_p(\la_1(p))$, where $\gamma_p$ is the curve in $T$ corresponding to $p\in M$.

Proceeding similarly to the previous case we get that in the limit $\delta\to 0$ equation \eqref{eqd2} in the new coordinates $r=r(q)$ reads
\begin{equation}\label{eqd3}
f_3f_{13}-f_1f_{33}+f_2f_{32}-f_3f_{22}-\frac{A}{2}f_2f_3=0.
\end{equation}
Substituting
$$
\left(\frac{A}{2}r_2\right)^2\mapsto e^{-Ar_2}
$$
equation \eqref{eqd3} takes the form
\begin{equation}\label{eqd3a}
f_3f_{13}-f_1f_{33}+e^{-Ar_2}(f_2f_{32}-f_3f_{22})=0.
\end{equation}
This is the most general equation in the family(C) of \cite{KP} (taking into account \cite[Remark 4.2]{KP}).

\section{Family (D)}
Consider a real-analytic hyper-CR Einstein-Weyl structure on $M=\R^3$. Then it prolongs to a hyper-CR Einstein-Weyl structure on $M_\C=\C^3$ (with a real structure). Let $T_\C$ be the corresponding complex twistor space \cite{H}. Then $T_\C$ fibers over $\C P^1=\C\cup\{\infty\}$, i.e. there is a projection $\pi\colon T_\C\to\C P^1$. Let us fix $\la_1,\la_2\in\C$ such that
$$
\la_1=a+\i b,\qquad \la_2=a-\i b,
$$
i.e. $\la_1=\bar\la_2$. Moreover, let $\la_3\in\R$ and $\la_4=\infty$. As in the real case we introduce new coordinates by formula
$$
z=(z_1,z_2,z_3)=(\gamma_z(\la_1),\gamma_z(\la_2),\gamma_z(\la_3)),
$$
where now $\gamma_z\colon\C P^1\to T_\C$ represents a point in $M_\C$ and the fibers $\pi^{-1}(\la_i)$, $i=1,2,3$, are identified with $\C$. Moreover, let 
$$
f(z)=\gamma_z(\la_4).
$$
Then $f$ is holomorphic and a curve $\gamma_z$ corresponds to a real point in $M\subset M_\C$ if
\begin{equation}\label{eqreal}
\gamma_z(\bar\la)=\bar\gamma_z(\la).
\end{equation}
In particular $f(z)\in\R$ if $z$ is a real point.

Let $\alpha_\la$ be the complex valued one-form defining null complex planes in $M_\C$. In the complex coordinates
$$
\alpha_\la=(\la-\la_2)(\la-\la_3)f_1dz_1+(\la-\la_1)(\la-\la_3)f_2dz_2+(\la-\la_1)(\la-\la_2)f_3dz_3,
$$
where now $f_i=\frac{\partial}{\partial z_i}f$. Denoting
$$
z_i=x_i+\i y_i
$$
we get from \eqref{eqreal} that on $M$
$$
x_1=x_2,\quad y_1=-y_2,\quad y_3=0.
$$
Writing
$$
f=g+\i h
$$
and exploiting the fact that $g_{x_i}=h_{y_i}$ and $g_{y_i}=-h_{x_i}$ we compute directly that on $M$
$$
\begin{aligned}
\alpha_\la=&(\la-c)(\la-a)((g_{x_1}+g_{x_2})(dx_1+dx_2)+(g_{y_1}-g_{y_2})(dy_1-dy_2))+\\
&(\la-c)b((g_{y_1}-g_{y_2})(dx_1+dx_2)+(g_{x_1}+g_{x_2})(dy_1-dy_2))+((\la-a)^2+b^2)g_{x_3}dx_3+\\
&\i(\la-c)(\la-a)(-(g_{y_1}+g_{y_2})(dx_1+dx_2)+(g_{x_1}-g_{x_2})(dy_1-dy_2))+\\
&\i(\la-c)b((g_{x_1}-g_{x_2})(dx_1+dx_2)+(g_{y_1}+g_{y_2})(dy_1-dy_2))-\i((\la-a)^2+b^2)g_{y_3}dx_3.
\end{aligned}
$$
Now, since $\partial_{x_1}+\partial_{x_2}$, $\partial_{y_1}-\partial_{y_2}$ and $\partial_{x_3}$ are tangent to $M$ and $h=0$ on $M$, the fact that $f$ is holomorphic implies that the imaginary part of $\alpha_\la$ on $M$ vanishes. Introducing 
$$
p_1=x_1+x_2,\quad p_2=y_1-y_2,\quad p_3=x_3
$$
we get that in the coordinates $p=(p_1,p_2,p_3)$ on $M$ the original Einstein-Weyl structure is defined by
\begin{equation}\label{formV33}
\alpha_\la=(\la-c)(\la-a)(g_1dp_1+g_2dp_2)+(\la-c)b(g_2dp_1+g_1dp_2)+((\la-a)^2+b^2)g_3dp_3,
\end{equation}
where $g_i=\frac{\partial}{\partial p_i}g$. Computing the integrability condition $d\alpha_\la\wedge\alpha_\la=0$ we get
\begin{equation}\label{eqD}
(a-c)(g_1g_{23}-g_2g_{13})+b(g_3g_{11}+g_3g_{22}-g_1g_{13}-g_2g_{23})=0
\end{equation}
which is the most symmetric equation in the family (D). To get a general equation in the family (D), with $a$ and $b$ depending on $p_1$ and $p_2$ and $c$ depending on $p_3$, one can take complex curves in $T_C$ instead of the fibers $\pi^{-1}(\lambda_1)$ and $\pi^{-1}(\lambda_2)$ and then proceed as in the real case.

\section{Hyper-Hermitian structures in dimension 4}
Let $T$ be the real twistor space corresponding to a hyper-Hermitian structure of split signature $(2,2)$ on a manifold $M$ (see \cite{C,DM,H,K}). Then $\dim T=3$ and, as in the Einstein-Weyl case, there is a natural fibration $\pi\colon T \to \R P^1$. Fixing two points $\lambda_1$ and $\lambda_2$ in $\R P^1$ we define coordinates of a point $p$
$$
p=(x,y)=(x^1,x^2,y^1,y^2)=(\gamma_p(\lambda_1),\gamma_p(\lambda_2))
$$ 
where $\gamma_p$ is the twistor curve in $T$ corresponding to $p$ and, similarly to the case of the Hirota equation, we assume that the fibers $\pi^{-1}(\lambda_i)$, $i=1,2$ are identified with $\R^2$. Choosing a third point $\lambda_3\in \R P^1$ and defining $f\colon M\to\R^2$ by formula
$$
f(p)=(f^1(p),f^2(p))=\gamma_p(\lambda_3)\in\R^2
$$
we get the following equation as the integrability condition (see \cite{K} for derivation)
\begin{equation}\label{eqH1}
\begin{aligned}
&f^i_{x_1y_1}(f^1_{x_2}f^2_{y_2}-f^2_{x_2}f^1_{y_2}) +f^i_{x_1y_2}(f^1_{x_2}f^2_{y_1}-f^2_{x_2}f^1_{y_1})\\
&+f^i_{x_2y_1}(f^1_{x_1}f^2_{y_2}-f^2_{x_1}f^1_{y_2}) +f^i_{x_2y_2}(f^1_{x_1}f^2_{y_1}-f^2_{x_1}f^1_{y_1})=0,\qquad i=1,2.
\end{aligned}
\end{equation}
Note that \eqref{eqH1} does not depend on the choice of $\lambda_i$, $i=1,2,3$. Now, assuming $\delta=\lambda_2-\lambda_1$ and substituting
$$
z^i=\frac{y^i-x^i}{\delta}
$$
in the limit $\delta\to 0$ we get
\begin{equation}\label{eqH2}
\begin{aligned}
&f^i_{z_1z_1}(f^1_{x_2}f^2_{z_2}-f^2_{x_2}f^1_{z_2}) +f^i_{z_1z_2}(f^1_{x_2}f^2_{z_1}-f^2_{x_2}f^1_{z_1}+f^1_{x_1}f^2_{z_2}-f^2_{x_1}f^1_{z_2}) +f^i_{z_2z_2}(f^1_{x_1}f^2_{z_1}-f^2_{x_1}f^1_{z_1})\\
&=f^i_{x_1z_2}(f^2_{z_2}f^1_{z_1}-f^1_{z_2}f^2_{z_1}) +f^i_{x_2z_1}(f^2_{z_1}f^1_{z_2}-f^1_{z_1}f^2_{z_2}),\qquad i=1,2.
\end{aligned}
\end{equation}
If additionally we assume that
$$
R^i=\frac{f^i-2\delta z^i-x^i}{\delta^2}
$$
in the limit $\delta\to 0$ we get
\begin{equation}\label{eqH3}
R^i_{z_1z_1}R^1_{z_2}+R^i_{z_2z_2}R^2_{z_1}-R^i_{z_1z_2}(R^2_{z_2}-R^1_{z_1})=2(R^i_{x_2z_1}-R^i_{x_1z_2}), \qquad i=1,2.
\end{equation}
This system coincides with (3.1) from \cite{K}.

\section{Veronese webs in dimension 4}
It is proved in \cite{K1} that the Veronese webs in dimension four locally correspond to the solutions to the system
\begin{equation}\label{sys1}
\begin{aligned}
&(\la_2-\la_1)f_0f_{12}+(\la_0-\la_2)f_1f_{02}+(\la_1-\la_0)f_2f_{01}=0,\\
&(\la_3-\la_1)f_0f_{13}+(\la_0-\la_3)f_1f_{03}+(\la_1-\la_0)f_3f_{01}=0,\\
&(\la_3-\la_2)f_0f_{23}+(\la_0-\la_3)f_2f_{03}+(\la_2-\la_0)f_3f_{02}=0,\\
\end{aligned}
\end{equation}
which is a generalization of the Hirota equation \eqref{eq1}. The system comes as the integrability condition
$$
d\alpha_\la\wedge\alpha=0
$$
with
$$
\alpha_\la=(\la-\la_0)(\la-\la_1)(\la-\la_2)(\la-\la_3)\sum_{i=0}^3\frac{(\la_4-\la_i)}{(\la-\la_i)}f_idp_i.
$$
The integrability condition means that the distributions $\ker\alpha_\la$ are integrable for each value of $\la$ and thus define a family of foliations which is exactly the Veronese web.

The twistor space $T$ for Veronese webs in any dimension is of dimension 2 and fibers over $\R P^1$ as in the three-dimensional case. The coordinates $(p_0,p_1,p_2,p_3)$ used above in \eqref{sys1} are defined, analogously to the three-dimensional case, by formula
$$
p_i=\gamma_p(\lambda_i),
$$
where $\gamma_p$ is the curve in $T$ corresponding to $p$ and $i=0,1,2,3$. Function $f$ is also defined as in the three-dimensional case by
$$
f(p)=\gamma_p(\lambda_4).
$$

Now, one can deform \eqref{sys1} first by assuming that some of $\la_i$ coincide and then replace the fibers $\pi^{-1}(\la_i)$ by arbitrary curves in $T$. We shall consider only two deformations which seem to be most interesting.

\paragraph{Case 1: $\la_i$, $i=0,1,2,3$, coincide.} Assume $\la_4=\infty$ and
$$
\la_i=\la_0+i\delta
$$
for $i=1,2,3$. We consider the following change of coordinates
$$
q_0=p_0,\quad q_1=\frac{p_1-p_0}{\delta},\quad q_2=\frac{p_2-2p_1+p_0}{\delta^2},\quad q_3=\frac{p_3-3p_2+3p_1-p_0}{\delta^3}.
$$
The transformation applied to all equations in the system \eqref{sys1} gives in the limit $\delta\to 0$ the following equation (up to a sign which is negative for the second equation)
$$
3(f_1f_{33}-f_3f_{13})-2(f_2f_{23}-f_3f_{22})=0.
$$
Then considering a sum of the first and the second equation (or the second and the third equation) we get in the limit $\delta\to 0$ the equation (up to a sign again which is negative for the sum of the second and the third equation)
$$
3(f_0f_{33}-f_3f_{03})-(f_2f_{13}-f_3f_{11})=0.
$$
Finally the sum of the first, the second multiplied by 2 and the third equation gives in the limit $\delta\to 0$ the equation
$$
2(f_0f_{23}-f_3f_{02})-(f_1f_{13}-f_3f_{11})=0.
$$
After rescaling $(q_0,q_1,q_2,q_3)\mapsto(q_0,q_1,2q_2,6q_3)$ we get the system
\begin{equation}\label{sys2}
\begin{aligned}
&f_1f_{33}-f_3f_{13}-f_2f_{23}+f_3f_{22}=0,\\
&f_0f_{33}-f_3f_{03}-f_2f_{13}+f_3f_{11}=0,\\
&f_0f_{23}-f_3f_{02}-f_1f_{13}-f_3f_{11}=0.
\end{aligned}
\end{equation}
The system is a direct generalization of the most symmetric equation in the family (C). One can check that it is equivalent to the integrability condition for
\begin{equation}\label{formV41}
\alpha_\la=(\la-\la_0)^3df-(\la-\la_0)^2(f_3dq_2+f_2dq_1+f_1dq_0)+(\la-\la_0)(f_3dq_1+f_2dq_0),
\end{equation}
which is the one-form defining the original Veronese web in the new coordinates. Note that the new form can be also derived as a limit $\delta\to 0$ of the original one-form defining the Veronese web. We have proved the following
\begin{theorem}
Locally, any Veronese web in dimension 4 can be put in the form \eqref{formV41}, where $f$ satisfies \eqref{sys2}. Conversely, any solution to \eqref{sys2} give rise to a Veronese web on $\R^4$ given by \eqref{formV41}.
\end{theorem}

\paragraph{Case 2: $\la_i$, $i=0,1,2,3,4$, coincide.} Now we assume that
$$
\la_4=\la_0+4\delta
$$
and additionally we consider
$$
H=\frac{f-4p_3+6p_2-4p_1+p_0}{\delta^4}.
$$
Then, proceeding as in the previous case and using the coordinates $q_i$, $i=0,1,2,3$ we get in the limit $\delta\to 0$ and after transformation $(q_0,q_1,q_2,q_3,H)\mapsto(q_0,q_1,2q_2,6q_3, 24H)$ the following system
\begin{equation}\label{sys3}
\begin{aligned}
&H_{13}-H_{22}+H_2H_{33}-H_3H_{23}=0,\\
&H_{03}-H_{12}+H_1H_{33}-H_3H_{13}=0,\\
&H_{02}-H_{11}+H_1H_{23}-H_2H_{13}=0.
\end{aligned}
\end{equation}
The system is a direct generalization of the hyper-CR equation \eqref{eq4}. One can check that it is equivalent to the integrability condition for the $\la$-dependent distribution
\begin{equation}\label{formV42}
D_\la=\spn\{\partial_i-\la(\partial_{i-1}+H_i\partial_3)\ |\ i=1,2,3\}.
\end{equation}
It is the tangent distribution of the corresponding Veronese web. We have proved the following
\begin{theorem}
Locally, any Veronese web in dimension 4 can be put in the form \eqref{formV42}, where $f$ satisfies \eqref{sys3}. Conversely, any solution to \eqref{sys3} give rise to a Veronese web on $\R^4$ given by \eqref{formV42}.
\end{theorem}

\paragraph{Remark.}
The two systems \eqref{sys2} and \eqref{sys3} are defined in the same coordinate system. Similarly to the 3-dimensional case their solutions are related by the following system
\begin{equation}\label{sys4}
\begin{aligned}
&f_0+H_1f_3=0,\\
&f_1+H_2f_3=0,\\
&f_2+H_3f_3=0.
\end{aligned}
\end{equation}
Namely, assuming that $f$ or $H$ is given then the consistency condition gives \eqref{sys2} or \eqref{sys3}, respectively.
\paragraph{Remark.}
All results in this section can be directly generalized to higher dimensions. In particular, the hierarchy \cite[eq. (6)]{DK} that extends the Hirota equation and system \eqref{sys1} above can be deformed (in the sense of the present section) to two other systems that appear to be the consistency conditions for the obvious extension of \eqref{sys4} given in the Introduction by formula \eqref{sys0}.

\end{document}